\documentclass{amsart}
\usepackage[matrix,arrow,cmtip]{xy}
\SelectTips{cm}{}
\usepackage{mathrsfs}
\usepackage{amssymb,amsthm,amsmath}

\numberwithin{equation}{section}

\newtheorem{Theorem}{Theorem}[section]
\newtheorem{Lemma}[Theorem]{Lemma}
\newtheorem{Proposition}[Theorem]{Proposition}

\theoremstyle{remark}
\newtheorem{Remark}[Theorem]{Remark}

\DeclareMathOperator{\Hom}{Hom}
\DeclareMathOperator{\Ext}{Ext}
\DeclareMathOperator{\End}{End}
\DeclareMathOperator{\Cok}{Cok}
\DeclareMathOperator{\Ker}{Ker}
\DeclareMathOperator{\Stab}{Stab}

\DeclareMathOperator{\Imm}{Im}

\newcommand{\Id}{\mathrm{Id}}
\newcommand{\id}{\mathrm{id}}
\newcommand{\Z}{\mathbb{Z}}
\newcommand{\C}{\mathbb{C}}
\newcommand{\Q}{\mathbb{Q}}
\newcommand{\R}{\mathbb{R}}
\newcommand{\GL}{\mathop\mathrm{GL}\nolimits}

\title{First extension groups of Verma modules and $R$-polynomials}
\author{Noriyuki Abe}
\address{Noriyuki Abe\\
Graduate School of Mathematical Sciences, the University of Tokyo,
 3--8--1 Komaba, Meguro-ku, Tokyo 153--8914, Japan.\\
abenori@ms.u-tokyo.ac.jp
}

\begin{document}
\maketitle

\begin{abstract}
We study the first extension groups between Verma modules.
There was a conjecture which claims that the dimensions of the higher extension groups between Verma modules are the coefficients of $R$-polynomials defined by Kazhdan-Lusztig.
This conjecture was known as the Gabber-Joseph conjecture (although Gabber and Joseph did not state.)
However, Boe gives a counterexample to this conjecture.
In this paper, we study how far are the dimensions of extension groups from the coefficients of $R$-polynomials.
\end{abstract}

\section{Introduction}
The category $\mathcal{O}$ is introduced by Bernstein-Gelfand-Gelfand and plays an important role in the representation theory.
One of the most important objects in $\mathcal{O}$ are the Verma modules and they are deeply investigated.

In this paper, we consider the $\Ext^i$-groups between Verma modules.
If $i = 0$, Verma~\cite{MR0218417} and Bernstein-Gelfand-Gelfand~\cite{MR0291204} determine the dimension of this group.
There are many studies about higher extension groups.
One of these studies is a work of Gabber-Joseph.
They proved some inequality between the dimensions of extension groups.
After that, it is conjectured that this inequality is, in fact, an equality.
Although not actually stated in \cite{MR644519}, this conjecture is known as the Gabber-Joseph conjecture.
If this conjecture is true, then the dimension of extension groups are the coefficients of $R$-polynomial.
However, Boe gives a counterexample to this conjecture~\cite{MR1197826}.

This conjecture is false even in the case of $i = 1$.
In this paper, we consider how far the dimensions of extension groups from the coefficients of $R$-polynomials.
Mazorchuk gives a formula of the dimension of the first extension group between Verma modules in a special case~\cite[Theorem~32]{MR2366357}.
Our formula is a generalization of his formula.

Now we state our main theorem.
Let $\mathfrak{g}$ be a semisimple Lie algebra over an algebraic closed field $K$ of characteristic zero.
Fix its Borel subalgebra $\mathfrak{b}$ and a Cartan subalgebra $\mathfrak{h}$.
Let $\Delta$ be the root system and $\rho$ the half sum of positive roots.
Fix a dominant integral element $\lambda\in \mathfrak{h}^*$. (It is sufficient to consider the integral case by a result of Soergel~\cite[Theorem~11]{MR1029692}.)
Let $M(x\lambda)$ be the Verma module with highest weight $x\lambda - \rho$ for $x\in W$.
Then by a result of Verma~\cite{MR0218417} and Bernstein-Gelfand-Gelfand~\cite{MR0291204}, if $x\ge y$, then there exists the unique (up to nonzero constant multiple) injective homomorphism $M(x\lambda)\to M(y\lambda)$.
Hence we can regard $M(x\lambda)$ as a submodule of $M(\lambda)$.
Then $M(w_0\lambda)$ is a submodule of $M(x\lambda)$ for all $x \in W$ where $w_0$ is the longest Weyl element.
Hence we have the homomorphism $\Ext^1(M(x\lambda),M(y\lambda))\to \Ext^1(M(w_0\lambda),M(\lambda))$.
Let $V_\lambda(x,y)$ be the image of this homomorphism.
We denote the unit element of $W$ by $e$.
Put $S_\lambda = \{s\in S\mid s(\lambda) = \lambda\}$.
For integral $\lambda,\mu\in\mathfrak{h}^*$, let $T_\lambda^\mu$ be the translation functor from $\lambda$ to $\mu$.
\begin{Theorem}[Proposition~\ref{prop:V is geometric repn}, Theorem~\ref{thm:calculation of V}, Theorem~\ref{thm:translation functor, surjective, kernel}]
Let $\lambda$ be a dominant integral element.
\begin{enumerate}
\item\label{enum:main theorem, has a W-module structure} If $\lambda$ is regular, then $V_\lambda(w_0,e)$ has a structure of $W$-module and it is isomorphic to $\mathfrak{h}^*$.
For $s\in S$, we denote the element in $V_\lambda(w_0,e)$ corresponding to the simple root whose reflection is $s$ by $v_s\in V_\lambda(w_0,e)$.
\item\label{enum:main theorem, determination of V_lambda, regular case} Assume that $\lambda$ is regular. For $x,y\in W$ and $s\in S$ such that $xs > x$, we have the following formula.
\begin{enumerate}
\item If $ys < y$, then $V_\lambda(xs,y) = s(V_\lambda(x,ys))$.
\item If $ys > y$, then $V_\lambda(xs,y) = Kv_s + s(V_\lambda(x,y))$.
\end{enumerate}
\item\label{enum:main theorem, determination of V_lambda, general} For general $\lambda$, the translation functor $T_\rho^\lambda$ induces a linear map $V_\rho(w_0,e)\to V_\lambda(w_0,e)$.
The kernel of this linear map is $\sum_{s\in S_\lambda}Kv_s$ and $T_\rho^\lambda(V_\rho(x,y)) = V_\lambda(x,y)$.
\end{enumerate}
\end{Theorem}
Notice that if $\lambda$ is regular and $x\not\ngeq y$, then $\Ext^1(M(x\lambda),M(y\lambda)) = 0$.
Hence $V_\lambda(x,y) = 0$.
Since $V_\lambda(x,x) = 0$ for all $x\in W$, we can determine the space $V_\lambda(x,y)$ inductively via \ref{enum:main theorem, determination of V_lambda, regular case} if $\lambda$ is regular.
Combining \ref{enum:main theorem, determination of V_lambda, general}, we can determine the space $V_\lambda(x,y)$.
In particular, we can calculate the dimension of $V_\lambda(x,y)$.
Namely, we can get the following theorem.
For $x,y\in W$, let $R_{x,y}$ be the polynomial defined in \cite[\S 2.]{MR560412}.
Let $\ell(x)$ be the length of $x\in W$.

\begin{Theorem}[Theorem~\ref{thm:Gabber-Joseph conjecture}]\label{thm:Main Theorem:Gabber-Joseph conjecture}
Assume that $\lambda$ is regular.
Then $\dim V_\lambda(x,y)$ is the coefficient of $q$ in $(-1)^{\ell(y) - \ell(x) - 1}R_{y,x}(q)$.
\end{Theorem}
In other words, $V_\lambda(x,y)$ satisfies the Gabber-Joseph conjecture.
Notice that the homomorphism $\Ext^1(M(x\lambda),M(y\lambda))\to \Ext^1(M(w_0\lambda),M(\lambda))$ is not injective.
It is easy to see that the kernel is isomorphic to $\Hom(M(x\lambda),M(\lambda)/M(y\lambda))$ if $x\ge y$ (Lemma~\ref{lem:preliminary of ext1}).

We summarize the contents of this paper.
In Section~\ref{sec:Preliminaries}, we gather preliminaries and prove some easy facts.
In Section~\ref{sec:Weyl group action}, we define a $W$-module structure on $\Ext^1(M(w_0\lambda),M(\lambda))$ for a dominant integral regular $\lambda$.
We also prove this module is isomorphic to $\mathfrak{h}^*$.
The proof of the main theorem in the regular case is done in Section~\ref{sec:Regular case}.
We finish the proof of the main theorem in Section~\ref{sec:Singular case}.

\section{Preliminaries}\label{sec:Preliminaries}
Let $\mathfrak{g}$ be a semisimple Lie algebra over an algebraic closed field $K$ of characteristic zero.
Fix a Borel subalgebra $\mathfrak{b}$ and a Cartan subalgebra $\mathfrak{h}\subset\mathfrak{b}$.
These determine the BGG category $\mathcal{O}$~\cite{MR0407097}.
We denote $\Hom_\mathcal{O}$ (resp.~$\Ext^i_\mathcal{O}$) by $\Hom$ (resp.~$\Ext^i$).
Denote the Weyl group of $\mathfrak{g}$ by $W$ and let $S$ be the set of its simple reflections.
For $w\in W$, let $\ell(w)$ be its length.
For $\lambda\in\mathfrak{h}^*$, let $\mathcal{O}_\lambda$ be the full-subcategory of $\mathcal{O}$ consisting of objects which have a generalized infinitesimal character $\lambda$.

In the rest of this paper, we only consider the objects which has a integral generalized infinitesimal character.
By~\cite[Theorem~11]{MR1029692}, it is sufficient to consider only the integral case.

For $\lambda\in\mathfrak{h}^*$, let $M(\lambda)$ be the Verma module with highest weight $\lambda - \rho$ where $\rho$ is the half sum of positive roots.
Assume that $\lambda$ is dominant integral.
Set $S_\lambda = S\cap\Stab_W\lambda$.
Then $W_{S_\lambda} = \Stab_W\lambda$ is generated by $S_\lambda$.
Put $W(S_\lambda) = \{w\in W\mid \text{$ws > w$ for all $s\in S_\lambda$}\}$.
This gives a complete representative set of $W/W_{S_\lambda}$.
For $x,y\in W(S_\lambda)$, $\Hom(M(x\lambda),M(y\lambda))\ne 0$ if and only if $x\ge y$.
Moreover, if $x\ge y$, then $\Hom(M(x\lambda),M(y\lambda))$ is one-dimensional and the nonzero homomorphisms from $M(x\lambda)$ to $M(y\lambda)$ are injective.
Hence we can regard $M(x\lambda)$ as a submodule of $M(y\lambda)$.
In the rest of this paper, we regard $M(x\lambda)$ as a submodule of $M(\lambda)$.
Then $M(w_0\lambda)$ is a submodule of $M(x\lambda)$ for all $x\in W$.

Let $\lambda,\mu\in\mathfrak{h}^*$ be integral dominant elements.
Then the translation functor $T_{\lambda}^\mu\colon \mathcal{O}_\lambda \to \mathcal{O}_\mu$ is defined~\cite{MR552943}.

Assume that $\lambda$ is a regular integral dominant element.
Then for $s\in S$, the wall-crossing functor $\theta^\lambda_s$ is defined by $\theta^\lambda_s = T_{\mu}^\lambda T_\lambda^\mu$ where $\Stab_W\mu = \{e,s\}$.
It is known that this functor is independent of $\mu$.
By the definition of $\theta^\lambda_s$, there are natural transformations $\Id\to \theta^\lambda_s$ and $\theta^\lambda_s\to \Id$.
Put $C^\lambda_s = \Cok(\Id\to \theta^\lambda_s)$ and $K_s^\lambda = \Ker(\theta_s^\lambda\to \Id)$.
Then $C^\lambda_s$ (resp.~$K_s^\lambda$) gives a right (resp.~left) exact functor from $\mathcal{O}_\lambda$ to $\mathcal{O}_\lambda$.
By the self-adjointness of $\theta^\lambda_s$, $(C^\lambda_s,K^\lambda_s)$ is an adjoint pair.
Moreover, the derived functor $LC_s$ gives an auto-equivalence of the derived category $D^b(\mathcal{O}_\lambda)$ and its quasi-inverse is $RK_s$.
For the action of $C_s,K_s$ on Verma modules, the following formulas hold:
Let $x\in W$ and $s\in S$ such that $xs > x$.
Then $C_s(M(x)) = M(xs)$, $K_s(M(xs)) = M(x)$ and there exists an exact sequence $0\to M(x)/M(xs)\to C_s(M(xs))\to M(xs)\to 0$.
Moreover, we have $L^kC_s(M(x)) = L^kC_s(M(xs)) = R^kK_s(M(xs)) = 0$ for all $k\ge 1$.

Let $\lambda$ be the dominant integral element of $\mathfrak{h}^*$.
Set
\[
	V_\lambda(x,y) = \Imm(\Ext^1(M(x\lambda),M(y\lambda))\to \Ext^1(M(w_0\lambda),M(\lambda)))
\]
for $x,y\in W$.
Put
\[
	V_\lambda = V_\lambda(w_0,e).
\]
\ref{enum:by big proj cover} of the following lemma is a part of the argument of the proof of~\cite[Lemma~33]{MR2366357}.
\ref{enum:dim of V_lambda} and \ref{enum:kernel of the map} follows from \ref{enum:by big proj cover}.
\begin{Lemma}\label{lem:preliminary of ext1}
Let $\lambda$ be a dominant integral element.
\begin{enumerate}
\item\label{enum:by big proj cover} Let $P$ be the projective cover of $M(w_0\lambda)$.
Let $M_i\subset P$ be the filtration such that $M_i/M_{i - 1} \simeq \bigoplus_{\ell(w) = i,w\in W(S_\lambda)}M(w\lambda)$.
Then we have
\begin{multline*}
\Ext^1(M(x\lambda),M(\lambda))\xleftarrow{\sim} \Hom(M(x\lambda),P/M_0)\\
\xleftarrow{\sim}\Hom(M(x\lambda),M_1/M_0)\simeq\bigoplus_{s\in S\setminus S_\lambda}\Hom(M(x\lambda),M(s\lambda)).
\end{multline*}
\item\label{enum:dim of V_lambda} We have $\dim V_\lambda = \#(S\setminus S_\lambda)$.\label{lem:preliminary of ext1:dim of V}
\item\label{enum:kernel of the map} For $x,y\in W$ such that $x\ge y$, the kernel of $\Ext^1(M(x\lambda),M(y\lambda))\to \Ext^1(M(w_0\lambda),M(\lambda))$ is isomorphic to $\Hom(M(x\lambda),M(\lambda)/M(y\lambda))$.
\end{enumerate}
\end{Lemma}
\begin{proof}
\ref{enum:by big proj cover} follows from the proof of~\cite[Lemma~33]{MR2366357}.
\ref{enum:dim of V_lambda} follows from \ref{enum:by big proj cover}.
We prove \ref{enum:kernel of the map}.
We have the long exact sequence
\begin{multline*}
0\to \Hom(M(x\lambda),M(y\lambda))\to \Hom(M(x\lambda),M(\lambda))
\\\to \Hom(M(x\lambda),M(\lambda)/M(y\lambda))\to \Ext^1(M(x\lambda),M(y\lambda))\to \Ext^1(M(x\lambda),M(\lambda))
\end{multline*}
The morphism $\Hom(M(x\lambda),M(y\lambda))\to \Hom(M(x\lambda),M(\lambda))$ is isomorphic by the classification of homomorphism between Verma modules.

Hence it suffices to prove that $\Ext^1(M(x\lambda),M(\lambda))\to \Ext^1(M(w_0\lambda),M(\lambda))$ is injective.
This follows from \ref{enum:by big proj cover}.
\end{proof}

\section{Weyl group action}\label{sec:Weyl group action}
In this section, fix a dominant regular integral element $\lambda$.
Put $V = V_\lambda$, $V(x,y) = V_\lambda(x,y)$, $\theta_s = \theta_s^\lambda$, $C_s = C_s^\lambda$ and $M(x) = M(x\lambda)$.
Since $M(x)\to \theta_s(M(x))\to M(x)$ is zero, we have $C_s(M(x))\to M(x)$.

We begin with the following lemma.
\begin{Lemma}\label{lem:existence of homs between ex seqs}
For $M_1,M_2,N_1,N_2$, take $a_i\in \Ext^1(M_i,N_i)$ for $i = 1,2$.
Take an exact sequence $0\to N_i\xrightarrow{k_i} X_i\xrightarrow{p_i} M_i\to 0$ corresponding to $a_i$.
Then for $f\colon N_1\to N_2$ and $g\colon M_1\to M_2$, there exists $\varphi$ such that the following diagram commutes if and only if $f_* a_1 = g^*a_2$:
\[
\xymatrix{
	0 \ar[r] & N_1\ar[r]\ar[d]^f & X_1\ar[r]\ar[d]^\varphi & M_1\ar[r]\ar[d]^g & 0\\
	0 \ar[r] & N_2\ar[r] & X_2\ar[r] & M_2\ar[r] & 0.
}
\]

\end{Lemma}
\begin{proof}
First assume that there exists $\varphi$.
Consider the following digram:
\[
\xymatrix{
0 \ar[r] & N_1\ar[r]\ar[d]^f & X_1\ar[r]\ar[d] & M_1\ar[r]\ar@{=}[d] &0\\
0 \ar[r] & N_2\ar[r] & X'\ar[r] & M_1\ar[r] &0,
}
\]
where the left square is a push-forward.
Then the exact sequence $0\to N_2\to X'\to M_1\to 0$ corresponds to $f_*a_1$ and the existence of $\varphi$ implies the existence of $\varphi'$ such that the following digram commutes:
\[
\xymatrix{
0 \ar[r] & N_2\ar[r]\ar@{=}[d] & X'\ar[r]\ar[d]^{\varphi'} & M_1\ar[r]\ar[d]^g &0\\
0 \ar[r] & N_2\ar[r] & X_2\ar[r] & M_2\ar[r] & 0.
}
\]
The same argument implies the existence of the following diagram:
\[
\xymatrix{
0 \ar[r] & N_2\ar[r]\ar@{=}[d] & X'\ar[r]\ar[d]^{\varphi''} & M_1\ar[r]\ar@{=}[d] &0\\
0 \ar[r] & N_2\ar[r] & X''\ar[r] & M_1\ar[r] & 0,
}
\]
where $0\to N_2\to X''\to M_1\to 0$ corresponds to $g^*a_2$.
The morphism $\varphi''$ should be isomorphism by $5$-Lemma.
Hence we have $f_*a_1 = g^*a_2$.

Conversely, we assume that $f_*a_1 = g^*a_2$.
Consider the following diagram:
\[
\xymatrix{
\Hom(M_1,N_2)\ar[r]^{k_{2*}}\ar[d]^{p_1^*} & \Hom(M_1,X_2) \ar[r]^{p_{2*}}\ar[d]^{p_1^*} & \Hom(M_1,M_2)\ar[r]^{\delta_2}\ar[d]^{p_1^*} & \Ext^1(M_1,N_2)\ar[d]^{p_1^*}\\
\Hom(X_1,N_2)\ar[r]^{k_{2*}}\ar[d]^{k_2^*} & \Hom(X_1,X_2) \ar[r]^{p_{2*}}\ar[d]^{k_2^*} & \Hom(X_1,M_2)\ar[r]^{\delta_1}\ar[d]^{k_2^*} & \Ext^1(X_1,N_2)\\
\Hom(N_1,N_2)\ar[r]^{k_{2*}}\ar[d]^{\delta_1} & \Hom(N_1,X_2) \ar[r]^{p_{2*}}\ar[d]^{\delta_1} & \Hom(N_1,M_2) & \\
\Ext^1(M_1,N_2)\ar[r]^{k_{2*}} & \Ext^1(M_1,X_2) & & .
}
\]
Recall that $f\in\Hom(N_1,N_2)$ and $g\in \Hom(M_1,M_2)$.
The assumption implies $\delta_1(f) = \delta_2(g)$.
Hence $\delta_1k_{2*}(f) = k_{2*}\delta_1(f) = k_{2*}\delta_2(g) = 0$.
Hence there exists an element $\varphi'\in \Hom(X_1,X_2)$ such that $k_2^*(\varphi') = k_{2*}(f)$.
Since $k_2^*p_{2*}(\varphi') = p_{2*}k_2^*(\varphi') = p_{2*}k_{2*}(f) = 0$, there exists $g'\in \Hom(M_1,M_2)$ such that $p_1^*(g') = p_{2*}(\varphi')$.
Then $p_{2*}(\varphi') = p_1^*(g')$ and $p_1^*(\varphi') = k_{2*}(f)$.
From the argument in the first part of this proof, we get $f_*a_1 = (g')^*a_2$, namely, we get $\delta_2(g') = \delta_1(f) = \delta_2(g)$.
Hence there exists $r\in \Hom(M_1,X_2)$ such that $p_{2*}(r) = g - g'$.
Set $\varphi = \varphi' + p_1^*(r)$.
Then $p_{2*}(\varphi) = p_1^*(g)$ and $p_1^*(\varphi) = k_{2*}(f)$.
This proves the lemma.
\end{proof}

The following lemma is well-known.
We give a proof for the sake of completeness.
\begin{Lemma}\label{lem:extension comming from translation}
Let $x\in W$ and $s\in S$ such that $x < xs$.
\begin{enumerate}
\item\label{enum:ext between sx and x} We have $\dim\Ext^1(M(xs),M(x)) = 1$.
The basis is given by the exact sequence $0\to M(x)\to\theta_s(M(x))\to M(xs)\to 0$.
\item\label{enum:map (sx,x)-->(w_0,e)} The homomorphism $\Ext^1(M(xs),M(x))\to \Ext^1(M(w_0),M(e))$ is injective and its image is independent of $x$.
\end{enumerate}
\end{Lemma}
\begin{proof}
Let $\mu$ be the integral dominant element such that $\Stab_W(\mu) = \{e,s\}$.
Then $\theta_s = T_\mu^\lambda T_\lambda^\mu$.

\ref{enum:ext between sx and x}
We have 
\[
\Ext^i(M(xs),\theta_s(M(x))) = \Ext^i(T_\lambda^\mu(M(xs)),T_\lambda^\mu(M(x)))
= \Ext^i(M(x\mu),M(x\mu))
\]
It is one-dimensional if $i = 0$ and zero if $i > 0$.
Hence from the exact sequence
\[
	0\to M(x)\to \theta_s(M(x))\to M(xs)\to 0,
\]
we have an exact sequence
\begin{multline*}
0\to \Hom(M(xs),M(x))\to \Hom(M(xs),\theta_s(M(x)))\\
\to \Hom(M(xs),M(xs)) \to \Ext^1(M(xs),M(x)) \to 0.
\end{multline*}
By the dimension counting, $\Hom(M(xs),M(x))\to \Hom(M(xs),\theta_s(M(x)))$ is isomorphic.
Hence $\Hom(M(xs),M(xs)) \simeq \Ext^1(M(xs),M(x))$.
The left hand side is one-dimensional.

\ref{enum:map (sx,x)-->(w_0,e)}
First we prove that $\Ext^1(M(xs),M(x))\to \Ext^1(M(w_0),M(e))$ is injective.
The kernel of the homomorphism is $\Hom(M(xs),M(e)/M(x))$ by Lemma~\ref{lem:preliminary of ext1} \ref{enum:kernel of the map}.
We have an exact sequence $0\to M(e)\to \theta_s(M(e))\to M(s)\to 0$ and $0\to M(x)\to \theta_s(M(x))\to M(xs)\to 0$.
Since $M(xs)\to M(s)$ is injective, we have that $M(e)/M(x)\to \theta_s(M(e)/M(x))$ is injective by the snake lemma.
Hence it is sufficient to prove that $\Hom(M(xs),\theta_s(M(e)/M(x))) = 0$.
We prove $\Hom(M(x\mu),M(\mu)/M(x\mu)) = 0$.
This follows from the following exact sequence
\begin{multline*}
0\to \Hom(M(x\mu),M(x\mu))\to \Hom(M(x\mu),M(\mu))\\
\to \Hom(M(x\mu),M(\mu)/M(x\mu))\to \Ext^1(M(x\mu),M(x\mu)) = 0.
\end{multline*}
We prove the independence of $x$.
We have the following diagram
\[
\xymatrix{
	M(xs)\ar[r]\ar[d] & M(x)\ar[d]\\
	M(s)\ar[r] & M(e).
}
\]
Hence we get the following diagram
\[
\xymatrix{
	0\ar[r] & M(x)\ar[r]\ar[d] & \theta_s(M(x))\ar[r]\ar[d] & M(xs)\ar[r]\ar[d] & 0\\
	0\ar[r] & M(e)\ar[r] & \theta_s(M(e))\ar[r] & M(s)\ar[r] & 0.
}
\]
Therefore the image of 
\[
	\Ext^1(M(xs),M(x))\to \Ext^1(M(xs),M(e))
\]
and
\[
	\Ext^1(M(s),M(e))\to \Ext^1(M(xs),M(e))
\]
coincide with each other by Lemma~\ref{lem:existence of homs between ex seqs}.
Therefore, we get \ref{enum:map (sx,x)-->(w_0,e)}
\end{proof}

Take a basis $v_s$ of $V(xs,x)$.

\begin{Lemma}
The set $\{v_s\mid s\in S\}$ is a basis of $V$.
\end{Lemma}
\begin{proof}
By Lemma~\ref{lem:preliminary of ext1} \ref{enum:by big proj cover}, we have the following commutative diagram
\[
\xymatrix{
\Ext^1(M(s),M(e)) \ar@{-}[r]^(0.44){\sim}\ar[d] & \bigoplus_{s'\in S}\Hom(M(s),M(s'))\ar[d]\\
\Ext^1(M(w_0),M(e))\ar@{-}[r]^(0.44){\sim} & \bigoplus_{s'\in S}\Hom(M(w_0),M(s')).
}
\]
We get the lemma.
\end{proof}

We need the following lemma to define the action of $s\in S$ on $V$.

\begin{Lemma}
The linear map $\Ext^1(M(w_0),M(e))\to \Ext^1(C_s(M(w_0)),M(e))$ is an isomorphism.
\end{Lemma}
\begin{proof}
By an exact sequence $0\to M(w_0s)/M(w_0)\to C_s(M(w_0))\to M(w_0)\to 0$, we get
\begin{multline*}
\Hom(M(w_0s)/M(w_0),M(e))\to \Ext^1(M(w_0),M(e))\\
\to \Ext^1(C_s(M(w_0)),M(e))\to \Ext^1(M(w_0s)/M(w_0),M(e))
\end{multline*}
It is sufficient to prove that $\Ext^i(M(w_0s)/M(w_0),M(e)) = 0$ for $i = 0,1$.
We have the following exact sequence
\begin{multline*}
0\to \Hom(M(w_0s)/M(w_0),M(e))\to \Hom(M(w_0s),M(e))\\
\to \Hom(M(w_0),M(e))\to \Ext^1(M(w_0s)/M(w_0),M(e))\\
\to \Ext^1(M(w_0s),M(e))\to \Ext^1(M(w_0),M(e)).
\end{multline*}
The homomorphism $\Hom(M(w_0s),M(e))\to \Hom(M(w_0),M(e))$ is an isomorphism.
The kernel of $\Ext^1(M(w_0s),M(e))\to \Ext^1(M(w_0),M(e))$ is isomorphic to $\Hom(M(w_0s),M(e)/M(e)) = 0$ by Lemma~\ref{lem:preliminary of ext1} \ref{enum:kernel of the map}.
We get the lemma.
\end{proof}
Using this lemma, we consider the following homomorphism:
\begin{align*}
V = \Ext^1(M(w_0),M(e)) & \to \Ext^1(LC_s(M(w_0)),LC_s(M(e)))\\
& \simeq \Ext^1(C_s(M(w_0)),M(s))\\
& \to \Ext^1(C_s(M(w_0)),M(e))\\
& \simeq \Ext^1(M(w_0),M(e)) = V.
\end{align*}
Define an action of $s$ on $V$ by the above homomorphism.
In other words, for $v_i\in V$ and the corresponding exact sequences $0\to M(e)\to X_i\to M(w_0)\to 0$, we have $s(v_1) = v_2$ if and only if there exists the following commutative diagram
\[
\xymatrix{
	0\ar[r] & C_s(M(e))\ar[r]\ar[d] & C_s(X_1)\ar[r]\ar[d] & C_s(M(w_0))\ar[r]\ar[d] & 0\\
	0\ar[r] & M(e)\ar[r] & X_2\ar[r] & M(w_0)\ar[r] & 0.
}
\]
Here, $C_s(M(e))\to M(e)$ and $C_s(M(w_0))\to M(w_0)$ are canonical homomorphisms.

The aim of this section is to show that this gives a structure of a $W$-module.
Since $\{C_s\}$ satisfies the braid relations~\cite[Lemma~5.10]{MR2139933}, this action satisfies the braid relations.

Take an exact sequence $0\to M(e)\to X\to M(w_0)\to 0$ and consider the following digram:
\[
\xymatrix{
	0 \ar[r]& M(e)\ar[r]\ar[d] & X\ar[r]\ar[d] & M(w_0)\ar[r]\ar[d] & 0\\
	0 \ar[r]& \theta_sM(e)\ar[r]\ar[d] & \theta_sX\ar[r]\ar[d] & \theta_sM(w_0)\ar[r]\ar[d] & 0\\
	0 \ar[r]& M(e)\ar[r] & X\ar[r] & M(w_0)\ar[r] & 0,
}
\]
here the vertical maps are natural transformations.
The compositions of left and right vertical morphisms are zero.
Hence composition of the middle vertical morphisms factors through $X\to M(w_0)$ and $M(e)\to X$.
By this way, we get an element of $\Hom(M(w_0),M(e))$.
This gives a morphism 
\[
	\alpha_s\colon V = \Ext^1(M(w_0),M(e))\to \Hom(M(w_0),M(e)).
\]
Since we fix an inclusion $M(w_0)\hookrightarrow M(e)$, we regard $\alpha_s(v)\in K$ for $v\in V$.
So we get $\alpha_s\in V^*$.

\begin{Lemma}\label{lem:varhi(a) = 0 implies s(a) = a}
If $v\in \Ext^1(M(w_0),M(e))$ satisfies $\alpha_s(v) = 0$, then $s(v) = v$.
\end{Lemma}
\begin{proof}
Let $0\to M(e)\to X\to M(w_0)\to 0$ be the corresponding exact sequence.
The assumption means that $X\to \theta_s(X)\to X$ is zero.
Hence $\theta_s(X)\to X$ factors through $\theta_s(X)\to C_s(X)$.
Namely, we get the following commutative diagram.
\[
\xymatrix{
0\ar[r] & C_s(M(e)) \ar[r]\ar[d] & C_s(X)\ar[r]\ar[d] & C_s(M(w_0))\ar[r]\ar[d] & 0,\\
0\ar[r] & M(e) \ar[r] & X\ar[r] & M(w_0)\ar[r] & 0.
}
\]
This means that $s(v) = v$.
\end{proof}

\begin{Lemma}\label{lem:varphi_s(v_s) = 2}
We have $\alpha_s(v_s) = 2$ and $s(v_s) = -v_s$.
\end{Lemma}
\begin{proof}
We consider $v_s$ as an element of $V(w_0,w_0s)$.
So we consider all things in $\{M\in \mathcal{O}_\lambda\mid \text{$[M:L(x)] = 0$ for $x\ne w_0,w_0s$}\}$.
This category is equivalent to the regular integral block of the BGG category of $\mathfrak{g} = \mathfrak{sl}_2(K)$.
So we may assume that $\mathfrak{g} = \mathfrak{sl}_2(K)$.

Set $P_0 = M(e)$, $P_1 = \theta_s(M(e))$.
Then $P_0\oplus P_1$ is a projective generator of $\mathcal{O}_\lambda$.
Set $A = \End(P_0\oplus P_1)$.
Then $\mathcal{O}_\lambda$ is equivalent to the category of finitely generated right $A$-modules.

We have that $\dim\Hom(P_0,P_0) = \dim\Hom(P_0,P_1) = \dim\Hom(P_1,P_0) = 1$, $\dim\Hom(P_1,P_1) = 2$.
Let $f \colon P_0\to P_1$ and $g \colon P_1\to P_0$ be the natural transformations.
Then $\Hom(P_0,P_1) = Kf$, $\Hom(P_1,P_0) = Kg$ and $\Hom(P_1,P_1) = K\id + Kfg$.
Moreover, $gf = 0$.
Set $e_i\colon P\to P_i$ be the projection.
Then $A = Ke_0 + Ke_1 + Kf + Kg + Kfg$, here $f$ stands for $P_0\oplus P_1\to P_0\xrightarrow{f}P_1\to P_0\oplus P_1$. ($g$ and $fg$ are similar.)

First, we calculate $\alpha_s(v_s)$.
To calculate it, we consider the following diagram:
\[
\xymatrix{
0\ar[r] & P_0\ar[r]^i\ar[d]^a & P_1\ar[r]^p\ar[d]^b & P_1/P_0\ar[r]\ar[d]^c & 0\\
0\ar[r] & P_1\ar[r]^j\ar[d]^s & P_1^{\oplus 2}\ar[r]^q\ar[d]^t & P_1\ar[r]\ar[d]^u & 0\\
0\ar[r] & P_0\ar[r]^i & P_1\ar[r]^p & P_1/P_0\ar[r] & 0,
}
\]
here $j$ is an inclusion $x\mapsto (x,0)$ and $q$ is second projection.
Notice that we used $\theta_s(P_0) = \theta_s(P_1/P_0) = P_1$ and the second exact sequence splits.

We regard $P_0$ and $P_1$ as right $A$-modules.
Then $P_0 = Ke_0 + Kg$, $P_1 = Ke_1 + Kf + Kfg$.
We regard $P_0$ as a submodule of $P_1$, namely, $P_0 = Kf + Kfg$.
We fix a inclusion $M(w_0) = P_1/P_0\to M(e) = P_0$ by $[e_1]\mapsto fg$, here $[e_1]$ is an image of $e_1$.
Set $h = fg$.
Since $\End(P_0) = K\id + Kh$, $b$ is given by $b = (\alpha_1 + \beta_1h,\alpha_2 + \beta_2h)$ for some $\alpha_1,\alpha_2,\beta_1,\beta_2\in K$.
Since $a$ and $i$ are both natural transformation.
So $a = i = f$.
By $bi = ja$, we have $\alpha_1 = 1$.
Since $c$ is given by $[e_1]\mapsto fg$, we get $\alpha_2 = 0$ and $\beta_2 = 1$ by $cp = qb$.
So we have $b = (1 + \beta_1h,h)$.

Next, we consider $t$.
Define $\gamma_1,\gamma_2,\delta_1,\delta_2\in K$ by $t = (\gamma_1 + \delta_1h,\gamma_2 + \delta_2h)$.
The morphism $s$ is given by $e_1\mapsto fg$, $f\mapsto 0$ and $fg\mapsto 0$.
Hence $tj = is$ implies that $\gamma_1 = 0$ and $\delta_1 = 1$.
Since $uq = pt$, we get $\gamma_2 = 1$.
Therefore, $t = (h,1 + \delta_2 h)$.
Hence the composition $tb$ is given by $2h$.
The image of $e_1$ under $P_1\to \theta_s(P_1)\to P_1$ is $2h$.
Since $P_1/P_0$ is given by $[e_1]\mapsto h$, this means $\alpha_s(v_s) = 2$.

Set $t' = (h,-1)\colon \theta_s(P_1)\simeq P_1^{\oplus 2}\to P_1$.
Then we have $t'b = 0$ and $(-u)q = pt'$.
This means that there exists a following diagram:
\[
\xymatrix{
0\ar[r] & C_s(P_0)\ar[r]\ar[d]^{\widetilde{s}} & C_s(P_1)\ar[r]\ar[d]^{\widetilde{t'}} & C_s(P_1/P_0)\ar[r]\ar[d]^{\widetilde{-u}} & 0\\
0\ar[r] & P_0\ar[r] & P_1\ar[r] & P_0/P_1\ar[r] & 0.
}
\]
Here $\widetilde{s}$ (resp.~$\widetilde{t'}$, $\widetilde{-u}$) is the morphism induced by $s$ (resp.~$t'$, $-u$).
Hence we have $s(v_s) = -v_s$.
\end{proof}

\begin{Lemma}\label{lem:s on V is a reflection}
For $v\in V$, we have $s(v) = v - \alpha_s(v)v_s$.
\end{Lemma}
\begin{proof}
Since $\alpha_s(v_s)\ne 0$, we have $V = Kv_s + \Ker \alpha_s$.
So we may assume that $v\in\Ker\alpha_s$ or $v = v_s$.
The first one is Lemma~\ref{lem:varhi(a) = 0 implies s(a) = a} and the second one is Lemma~\ref{lem:varphi_s(v_s) = 2}.
\end{proof}

\begin{Proposition}
The action of $s\in S$ on $V$ defines a representation of $W$.
\end{Proposition}
\begin{proof}
We should prove $s^2 = 1$.
This follows from Lemma~\ref{lem:s on V is a reflection} and Lemma~\ref{lem:varphi_s(v_s) = 2}.
\end{proof}

In fact, $V$ is isomorphic to $\mathfrak{h}^*$ as a $W$-module.
We prove it in the rest of this section.

\begin{Lemma}\label{lem:s(V(x,y)), easy case}
If $xs > x$ and $ys > y$, then $V(xs,ys) = s(V(x,y))$.
\end{Lemma}
\begin{proof}
This follows from $C_s(M(x))\simeq M(xs)$ and $C_s(M(y))\simeq M(ys)$.
\end{proof}

\begin{Lemma}\label{lem:action of s is cos(pi/m)}
Let $m'_{s,s'}$ be an order of $ss'\in \GL(V)$.
Then there exists $v'_s\in K^\times v_s$ such that $s(v'_{s'}) = v'_{s'} - 2\cos(\pi/m'_{s,s'})v_s$ for all $s,s'\in S$.
\end{Lemma}
\begin{proof}
Since $W$ is a finite group, $V$ is defined over $\overline{\Q}$.
Hence we may assume $K = \overline{\Q}$.
Moreover, by the base change, we may assume that $K = \C$.
Since $W$ is a finite group, $V$ has a $W$-invariant inner product $\langle\cdot,\cdot\rangle$.
Take $v''_s\in \C v_s$ such that $\langle v''_s,v''_s\rangle = 2$.
By Lemma~\ref{lem:s on V is a reflection}, $s(v^*) = v^* - v^*(v_s)\alpha_s$ for $v^*\in V$.
Hence the kernel of $s + 1\colon V^*\to V^*$ is $K\alpha_s$.
The linear map $v\mapsto \langle v,v''_s\rangle$ belongs to $\Ker(s + 1) = K\alpha_s$.
Since $\langle v''_s,v''_s\rangle = \alpha_s(v_s)$, we have $s(v) = v - \langle v,v''_s\rangle v''_s$ by Lemma~\ref{lem:s on V is a reflection}.

Put $a_{s,s'} = \langle v''_s,v''_{s'}\rangle$.
Then we have $a_{s,s'} = \overline{a_{s',s}}$.
By Lemma~\ref{lem:s on V is a reflection}, $Kv'_s + Kv'_{s'}$ is stable under the action of $ss'$.
Its characteristic polynomial is $t^2 + (2 - \lvert a_{s,s'}\rvert^2)t + 1$.
Let $\alpha,\alpha^{-1}$ be its eigenvalue.
Since $(W,S)$ is a Weyl group, $m'_{s,s'} \le 6$.
Hence $\alpha = e^{2\pi\sqrt{-1}/m'_{s,s'}}$ or $\alpha = e^{-2\pi\sqrt{-1}/m'_{s,s'}}$.
We get $\lvert a_{s,s'}\rvert^2 = 2 + \alpha + \alpha^{-1} = 2 + 2\cos2\pi/m'_{s,s'} = (2\cos(\pi/m'_{s,s'}))^2$.
Hence we have $a_{s,s'} = 2e^{\sqrt{-1}\theta_{s,s'}}\cos(\pi/m'_{s,s'})$ for some $\theta_{s,s'}\in\R$.
Recall that the Coxeter graph of $(W,S)$ is a tree.
Hence we can choose $v'_s\in e^{\sqrt{-1}\R}v''_s$ such that $a_{s,s'} = 2\cos(\pi/m'_{s,s'})$.
So we get the lemma.
\end{proof}

\begin{Lemma}\label{lem:condition for v_s in Ext}
For $x,y\in W$ such that $xs > x$, $x \ge y$ and $y < ys$, let $A$ (resp.~$B$) be the image of $\Ext^1(M(xs),M(ys))\to \Ext^1(M(xs),M(y))$ (resp.~$\Ext^1(M(xs),M(x))\to \Ext^1(M(xs),M(y))$).
Then $A\cap B = 0$ if and only if $x\not\ge ys$.
\end{Lemma}
\begin{proof}
By Lemma~\ref{lem:extension comming from translation}, $B$ is one-dimensional.
If $x \ge ys$, then we have a homomorphism $\Ext^1(xs,x)\to \Ext^1(xs,ys)$.
So $B\subset A$.

On the other hand, assume that $A\cap B \ne 0$.
Let $\mu\in \mathfrak{h}^*$ be a dominant integral element such that $\Stab_W(\mu) = \{e,s\}$.
Then $T_\lambda^\mu$ induces an homomorphism $\Ext^1(M(xs),M(y))\to \Ext^1(M(xs\mu),M(y\mu))$.
Since $T_\lambda^\mu(M(xs)) = T_\lambda^\mu(M(x)) = M(x\mu)$, the image of $B$ under $T_\lambda^\mu$ is zero.
Hence the homomorphism 
\[
	\Ext^1(M(xs),M(ys))\to \Ext^1(M(xs\mu),M(y\mu))
\]
has a kernel.
This homomorphism is equal to 
\[
	\Ext^1(M(xs),M(ys))\to \Ext^1(M(xs),\theta_s(M(y))).
\]
From an exact sequence 
\[
	0\to M(ys)\to \theta_s(M(y))\to C_s(M(ys))\to 0, 
\]
we have $\Hom(M(xs),C_s(M(ys))\ne 0$.
By the adjointness, 
\[
	\Hom(M(x),C_s(M(ys))) = \Hom(K_s(M(xs)),M(ys)) = \Hom(M(x),M(ys)).
\]
Hence we have $x\ge ys$.
\end{proof}

\begin{Lemma}\label{lem:rank 2 case, injective}
For $s,s'\in S$, let $\langle s,s'\rangle$ be the group generated by $\{s,s'\}$.
Then for all $x,y\in w_0\langle s,s'\rangle$, $\Ext^1(M(x),M(y))\to \Ext^1(M(w_0),M(e))$ is injective.
\end{Lemma}
\begin{proof}
Denote the irreducible quotient of $M(x)$ by $L(x)$.
Let $w\in \langle s,s'\rangle$ be the longest element.
Put 
\[
	\mathcal{O}' = \{M\in \mathcal{O}_\lambda\mid \text{$[M:L(z)] = 0$ for $z\not\in w_0\langle s,s'\rangle$}\}.
\]
Then $\mathcal{O}'$ is equivalent to the regular integral block of the BGG category of semisimple Lie algebra of rank $2$.
Applying Lemma~\ref{lem:preliminary of ext1} \ref{enum:kernel of the map} to this category, the kernel of 
\[
	\Ext^1(M(x),M(y))\to \Ext^1(M(w_0),M(w_0w))
\]
is isomorphic to $\Hom(M(x),M(w_0w)/M(y))$.
However, we have that $[M(w_0w):L(x)] = [M(y):L(x)] = 1$.
Hence this space is zero.
By Lemma~\ref{lem:preliminary of ext1} \ref{enum:kernel of the map}, the kernel of 
\[
	\Ext^1(M(w_0),M(w_0w))\to \Ext^1(M(w_0),M(e))
\]
is isomorphic to $\Hom(M(w_0),M(e)/M(w_0w))$.
It is zero since $[M(e):L(w_0)] = [M(w_0w):L(w_0)] = 1$.
\end{proof}

\begin{Remark}\label{rem:when v_s in V}
By Lemma~\ref{lem:condition for v_s in Ext} and Lemma~\ref{lem:rank 2 case, injective}, we have the following.
For $x,y\in w_0\langle s,s'\rangle$ such that $xs > x \ge y < ys$, $v_s\in V(xs,ys)$ if and only if $x\ge ys$.
Since $s(v_s) = -v_s$ and $V(xs,ys) = s(V(x,y))$ (Lemma~\ref{lem:s(V(x,y)), easy case}), under the same conditions, we have $v_s\in V(x,y)$ if and only if $x\ge ys$.
\end{Remark}

\begin{Proposition}\label{prop:V is geometric repn}
The $W$-representation $V$ is isomorphic to the geometric representation defined in~\cite{MR1890629}.
\end{Proposition}
\begin{proof}
Let $w' \in \langle s,s'\rangle$ be the longest element and put $w = w_0w'$.
Let $m_{s,s'}$ be an order of $ss'\in W$.
It is sufficient to prove that $m_{s,s'}$ is an order of $ss'\in \GL(V)$.

Fix $s,s'\in S$ and set $m = m_{s,s'}$.
Put $V' = Kv_s + Kv_{s'} = V(w_0,w)$.
Then $V'$ is $\langle s,s'\rangle$-stable.
Let $n$ be an order of $ss'\in \GL(V')$.
Then $n|m$.
We prove $n = m$.
Since $(W,S)$ is a Weyl group, $m = 2,3,4,6$.
If $m = 2,3$, then there is nothing to prove.

Assume that $n = 2$ and $m\ne 2$.
Then by Lemma~\ref{lem:action of s is cos(pi/m)}, $s'v'_{s} = v'_{s}$.
Therefore, $v_s\in s'V(ws,w) = V(wss',ws')$ by Lemma~\ref{lem:s(V(x,y)), easy case}.
By Remark~\ref{rem:when v_s in V}, $ss'\ge s's$.
This is a contradiction since $m \ne 2$.

Assume that $(n,m) = (3,6)$.
Then we have $s's(v'_{s'}) = -v'_{s}$.
Hence, $v_s\in s'sV(ws',w) = V(ws'ss',wss')$.
By Remark~\ref{rem:when v_s in V}, $s'ss'\ge ss's$.
This is a contradiction since $m ~= 6$.
\end{proof}
This proposition says that $V\simeq \mathfrak{h}^*$ as $W$-module.
Let $\alpha$ be a simple root.
By the proof, $v_{s_\alpha}$ corresponds to $\alpha$ up to constant multiple.

\section{Regular case}\label{sec:Regular case}
As in the previous section, fix a regular integral dominant element $\lambda$.
We continue to use the notation in the previous section.
In this section, we determine $V(x,y)\subset V$.

We need the graded BGG category.
By a result of Beilinson-Ginzburg-Soergel~\cite{MR1322847}, there exists a graded algebra $A = \bigoplus_{i\ge 0}A_i$ such that $\mathcal{O}_\lambda$ is equivalent to the category of $A$-modules.
Let $\widetilde{\mathcal{O}_\lambda}$ be the category of graded $A$-modules.
Then we have the forgetful functor $\widetilde{\mathcal{O}_\lambda}\to \mathcal{O}_\lambda$.
For $M = \bigoplus_iM_i\in\widetilde{\mathcal{O}_\lambda}$, define $M\langle n\rangle\in\widetilde{\mathcal{O}_\lambda}$ by $(M\langle n\rangle)_i = M_{i - n}$.
There exists a module $\widetilde{M}(x)\in \widetilde{\mathcal{O}_\lambda}$ such that its image in $\mathcal{O}_\lambda$ is isomorphic to $M(x)$ and it is unique up to grading shift~\cite[3.11]{MR1322847}.
Take a grading of $\widetilde{M}(x)$ such that its top is in degree $0$.
Set $\widetilde{M}^0(x) = \widetilde{M}(x)\langle \ell(x)\rangle$.
Then by the proof of \cite[Corollary~23]{MR2366357}, if $x\ge y$, then there exists a degree zero injective homomorphism $\widetilde{M}^0(x)\to \widetilde{M}^0(y)$.

For $M,N\in \widetilde{\mathcal{O}_\lambda}$, let $\Hom(M,N)_0$ be the space of homomorphisms of degree zero and $\Ext^i(M,N)_0$ its derived functors.
The following lemma is proved in \cite[Theorem~32]{MR2366357}.
Howerver, the author thinks, although his proof is correct, the shift in the statement of the theorem is wrong.
So we give a proof (it is the same as in the proof in \cite{MR2366357}).

\begin{Lemma}\label{lem:Mazorchuk's calculation}
We have $\Ext^1(\widetilde{M}^0(x),\widetilde{M}^0(e)\langle i\rangle)_0 = 0$ if $i\ne 2$.
Moreover, if $i = 2$, its dimension is equal to $\#\{s\in S\mid s\le x\}$.
\end{Lemma}
\begin{proof}
We use the argument in the proof of Lemma~\ref{lem:preliminary of ext1} with a grading.
Let $\widetilde{P}$ be the projective cover of $\widetilde{M}(w_0)$.
Then it has a Verma flag $M_i$ such that $M_i/M_{i - 1}$ is isomorphic to $\bigoplus_{\ell(x) = i}M(x)$ (here, we ignore the grading).
The by Lemma~\ref{lem:preliminary of ext1} \ref{enum:by big proj cover}, we have $\Hom(\widetilde{M}^0(x\lambda),M_1/M_0)_0\simeq \Ext^1(\widetilde{M}^0(x\lambda),M_0)$.
We consider the grading of $M_0$ and $M_1/M_0$.
Take $k\in \Z$ such that $M_0 \simeq \widetilde{M}(e)\langle k\rangle$.
Then the multiplicity of $\widetilde{M}(w_0)\langle k\rangle$ in $\widetilde{M}(e)$ is $1$ by \cite[Theorem~3.11.4]{MR1322847}.
Since $\widetilde{M}(w_0)\langle \ell(w_0)\rangle\subset \widetilde{M}(e)$, we have $k = \ell(w_0)$.
By the same argument, we have $M_1/M_0\simeq \bigoplus_s\widetilde{M}(s)\langle \ell(w_0) - 1\rangle$.
So we have 
\[
	\Ext^1(\widetilde{M}^0(x),\widetilde{M}^0(e)\langle  i\rangle)
	\simeq \Hom\left(\widetilde{M}^0(x),\bigoplus_{s\in S}\widetilde{M}^0(s)\langle i - 2\rangle\right).
\]
We get the lemma.
\end{proof}

We use the following abbreviations.
\begin{itemize}
\item $E^i(x\langle k\rangle,y\langle l\rangle) = \Ext^i(\widetilde{M}^0(x)\langle k\rangle,\widetilde{M}^0(y)\langle l\rangle)_0$.
\item $E^i(x\langle k\rangle,y/z\langle l\rangle) = \Ext^i(\widetilde{M}^0(x)\langle k\rangle,\widetilde{M}^0(y)/\widetilde{M}^0(z)\langle l\rangle)_0$.
\end{itemize}
By the above lemma, $V(x,y)$ is the image of $E^1(x,y\langle 2\rangle)\to E^1(w_0,e\langle 2\rangle)$.

\begin{Lemma}\label{lem:vanish of ext1, by grading}
If $i\le 1$, then $\Ext(\widetilde{M}^0(x),\widetilde{M}^0(y)\langle i\rangle)_0 = 0$.
\end{Lemma}
\begin{proof}
We have $E^1(x,y\langle i\rangle) = E^0(x,e/y\langle i\rangle)$ by Lemma~\ref{lem:preliminary of ext1} \ref{enum:kernel of the map} and the above lemma.
Let $L$ be the unique irreducible quotient of $\widetilde{M}^0(x)$.
By \cite[Theorem~3.11.4]{MR1322847}, we have $[\widetilde{M}^0(e)/\widetilde{M}(y)\langle i\rangle : L] = 0$.
Hence we have $E^0(x,e/y\langle i\rangle) = 0$.
\end{proof}

There is a graded lift of $\theta_s$~\cite{MR2005290}.
We take a lift $\widetilde{\theta_s}$ of $\theta_s$ such that $\widetilde{\theta_s}$ is self-dual and there exist degree zero natural transformations $\Id\langle 1\rangle \to\widetilde{\theta_s}\to \Id\langle -1\rangle$.
Let $\widetilde{C_s}$ (resp.~$\widetilde{K_s}$) be the cokernel (resp.~kernel) of $\Id\langle 1\rangle \to \widetilde{\theta_s}$ (resp.~$\widetilde{\theta_s}\to \Id\langle -1\rangle$).
Then $(\widetilde{C_s},\widetilde{K_s})$ is an adjoint pair.
By \cite[Theorem~3.6]{MR2005290}, we have the following formulas for $x\in W$ and $s\in S$ such that $xs > x$:
\begin{itemize}
\item $\widetilde{C_s}(\widetilde{M}^0(x)) = \widetilde{M}^0(xs)\langle -1\rangle$ and $L^k\widetilde{C_s}(\widetilde{M}^0(x)) = 0$ for $k > 0$.
\item We have an exact sequence $0\to (\widetilde{M}^0(x)/\widetilde{M}^0(xs))\langle 1\rangle\to C_s(\widetilde{M}^0(xs))\to \widetilde{M}^0(xs)\langle -1\rangle \to 0$.
\item $\widetilde{K_s}(\widetilde{M}^0(xs)) = \widetilde{M}^0(x)\langle 1\rangle$ and $R^k\widetilde{K_s}(\widetilde{M}^0(xs)) = 0$ for $k > 0$.
\end{itemize}
We also have that $L\widetilde{C_s}$ gives an auto-equivalence of $D^b(\widetilde{\mathcal{O}_\lambda})$ and its quasi-inverse functor is $R\widetilde{K_s}$.

Now we prove the main theorem in the regular case.
We have already proved \ref{enum:calculation of V, first case} (Lemma~\ref{lem:s(V(x,y)), easy case}).
\begin{Theorem}\label{thm:calculation of V}
Let $x,y\in W$ and $s\in S$ such that $xs > x\ge y$.
\begin{enumerate}
\item\label{enum:calculation of V, first case} If $ys < y$, then $V(xs,y) = s(V(x,ys))$.
\item\label{enum:calculation of V, second case} If $ys > y$, then $V(xs,y) = Kv_s + s(V(x,y))$.
\end{enumerate}
\end{Theorem}
\begin{proof}
From \ref{enum:calculation of V, first case}, we have $s(V(x,y)) = V(xs,ys)$.
Since $y < ys$, we have a homomorphism $\Ext^1(M(xs),M(ys))\to \Ext^1(M(xs),M(y))$.
Hence $V(xs,ys)\subset V(xs,y)$.
We also have $\Ext^1(M(xs),M(x))\to \Ext^1(M(xs),M(y))$.
Since $Kv_s$ is the image of $\Ext^1(M(xs),M(x))$, the right hand side is contained in the left hand side.

From an exact sequence
\[
	0\to (\widetilde{M}^0(y)/\widetilde{M}^0(ys))\langle 2\rangle \to \widetilde{C_s}(\widetilde{M}^0(ys))\langle 1\rangle \to \widetilde{M}^0(ys)\to 0,
\]
we have an exact sequence
\begin{multline*}
\Ext^0(\widetilde{M}^0(xs),\widetilde{C_s}(\widetilde{M}^0(ys))\langle 1\rangle)_0
\to E^0(xs,ys)\\ \to E^1(xs,y/ys\langle 2\rangle)
\to \Ext^1(\widetilde{M}^0(xs),\widetilde{C_s}(\widetilde{M}^0(ys))\langle 1\rangle)_0.
\end{multline*}

Since $R\widetilde{K_s}$ is the quasi-inverse functor of $L\widetilde{C_s}$, for $i\ge 0$, we have
\begin{align*}
	\Ext^i(\widetilde{M}^0(xs),C_s(\widetilde{M}^0(ys))\langle 1\rangle)_0
	& = \Hom_{D^b(\widetilde{\mathcal{O}_\lambda})}(\widetilde{M}^0(xs),L\widetilde{C_s}(\widetilde{M}^0(ys))\langle 1\rangle[i])_0\\
	& = \Hom_{D^b(\widetilde{\mathcal{O}_\lambda})}(R\widetilde{K_s}(\widetilde{M}^0(xs)),\widetilde{M}^0(ys)\langle 1\rangle[i])_0\\
	& = \Hom_{D^b(\widetilde{\mathcal{O}_\lambda})}(\widetilde{M}^0(x)\langle 1\rangle,\widetilde{M}^0(ys)\langle 1\rangle[i])_0\\
	& = \Ext^i(\widetilde{M}^0(x),\widetilde{M}^0(ys))_0.
\end{align*}
If $i = 1$, then this is zero by Lemma~\ref{lem:vanish of ext1, by grading}.
Hence we get an exact sequence
\[
	E^0(x,ys)\to E^0(xs,ys)\to E^1(xs,y/ys\langle 2\rangle)\to 0.
\]

Assume that $x \not\ge ys$.
Then $E^0(x,ys) = 0$.
Therefore, we have that $\dim E^1(xs,y/ys\langle 2\rangle) = \dim E^0(xs,ys) = 1$.
From an exact sequence
\[
	E^1(xs,ys\langle 2\rangle)\to E^1(xs,y\langle 2\rangle)\to E^1(xs,y/ys\langle 2\rangle),
\]
the codimension of the image $A$ of $E^1(xs,ys\langle 2\rangle)\to E^1(xs,y\langle 2\rangle)$ is less than or equal to $1$.
We also have that the image $B$ of $E^1(xs,x\langle 2\rangle)\to E^1(xs,ys\langle 2\rangle)$ is one-dimensional by Lemma~\ref{lem:extension comming from translation}.
By Lemma~\ref{lem:condition for v_s in Ext}, $A\cap B = 0$.
Hence $A + B = E^1(xs,y\langle 2\rangle)$.
This implies the theorem in this case.

Next assume that $x\ge ys$.
Then $E^0(x,ys)\to E^0(xs,ys)$ is isomorphic.
Hence $E^1(xs,y/ys\langle 2\rangle) = 0$.
Therefore, $E^1(xs,ys\langle 2\rangle)\to E^1(xs,y\langle 2\rangle)$ is surjective.
This implies $V(xs,ys)\supset V(xs,y)$.
\end{proof}

As a corollary, we can determine the dimension of $V(x,y)$.
Let $R_{y,x}(q)$ be the polynomial defined in \cite{MR560412}.

\begin{Theorem}\label{thm:Gabber-Joseph conjecture}
Assume that $\lambda$ is regular.
Then $\dim V_\lambda(x,y)$ is the coefficient of $q$ in $(-1)^{\ell(y) - \ell(x) - 1}R_{y,x}(q)$.
\end{Theorem}
\begin{proof}
Put $n_{y,x} = \dim V(x,y)$.
Then by the above theorem and its proof, for $x,y\in W$, $s\in S$ such that $xs > x \ge y$, we have
\begin{itemize}
\item If $ys < y$, then $n_{xs,y} = n_{x,ys}$.
\item If $ys > y$ and $x\ge ys$, then $n_{xs,y} = n_{x,y}$.
\item If $ys > y$ and $x\not\ge ys$, then $n_{xs,y}\le n_{x,y} + 1$.
\end{itemize}
Let $r_{y,x}$ be the coefficient of $q$ in $(-1)^{\ell(y) - \ell(x) - 1}R_{y,x}(q)$.
By \cite[(2.0.b), (2.0.c)]{MR560412}, the constant term of $(-1)^{\ell(y) - \ell(x) - 1}R_{y,x}(q)$ is $0$ or $1$ and it is $1$ if and only if $y\le x$.
Hence using \cite[(2.0.b), (2.0.c)]{MR560412}, we have
\begin{itemize}
\item If $ys < y$, then $r_{xs,y} = r_{x,ys}$.
\item If $ys > y$ and $x\ge ys$, then $r_{xs,y} = r_{x,y}$.
\item If $ys > y$ and $x\not\ge ys$, then $r_{xs,y} = r_{x,y} + 1$.
\end{itemize}
Hence we get $n_{y,x}\le r_{y,x}$.
To prove $n_{y,x} = r_{y,x}$, it is sufficient to prove $n_{w_0,x} = r_{w_0,x}$.
We prove this by backward induction on $\ell(x)$.
By Lemma~\ref{lem:Mazorchuk's calculation}, we have $n_{w_0,x} = \#\{s'\in S\mid w_0s'\ge x\}$.

Take $s\in S$ such that $xs > x$.
If $w_0s \ge xs$, then $r_{w_0,x} = r_{w_0s,x} = r_{w_0,xs} = n_{w_0,xs}$.
If $w_0s\not\ge xs$, then $r_{w_0,x} = r_{w_0s,x} + 1 = r_{w_0,xs} + 1 = n_{w_0,xs} + 1$.
We compare $X = \{s'\in S\mid w_0s'\ge x\}$ and $Y = \{s'\in S\mid w_0s'\ge xs\}$.
Since $xs > x$, we have $Y\subset X$.

Assume $s'\in Y$ and $s'\ne s$.
Then $s' \le w_0xs$.
Hence $s'$ appears in a reduced expression of $w_0xs$.
Since $s\ne s'$, $s'$ appears in a reduced expression of $w_0x$.
Hence $s'\in X$.
This implies $X\cap (S\setminus\{s\}) = Y\cap (S\setminus \{s\})$.

Since $xs > x$, we have $w_0s\ge x$.
Therefore, $s\in X$.
Hence if $w_0s\ge xs$, we have $X = Y$, this implies $n_{w_0,x} = n_{w_0,xs}$.
If $w_0s\not\ge xs$, we have $X = Y\amalg \{s\}$, this implies $n_{w_0,x} = n_{w_0,xs} + 1$.
Therefore, $r_{w_0,x} = n_{w_0,x}$.
\end{proof}

\section{Singular case}\label{sec:Singular case}
In this section, we fix a dominant integral (may be singular) element $\lambda\in \mathfrak{h}^*$.
We also fix a regular integral dominant element $\lambda_0\in \mathfrak{h}^*$.
Then the translation functor $T_{\lambda_0}^\lambda$ is defined and it gives $V_{\lambda_0}(x,y)\to V_\lambda(x,y)$.
Recall the notation $S_\lambda = \{s\in S\mid s(\lambda) = \lambda\}$.

In this section, we prove the following theorem.
\begin{Theorem}\label{thm:translation functor, surjective, kernel}
\begin{enumerate}
\item\label{enum:map by translation functor is surjective} The homomorphism $V_{\lambda_0}(x,y)\to V_\lambda(x,y)$ induced by the translation functor is surjective.
\item\label{enum:kernel of the map by translation functor} The kernel of $V_{\lambda_0}\to V_\lambda$ is $\sum_{s\in S_\lambda}Kv_s$.
\end{enumerate}
\end{Theorem}

We use the notation $\widetilde{\mathcal{O}_\lambda}$, $M\langle n\rangle$ and $\Hom(M,N)_0$ which we use in the previous section.
Then using the same argument in \cite{MR2005290}, $T_{\lambda_0}^\lambda$ and $T_\lambda^{\lambda_0}$ have graded lifts $\widetilde{T_{\lambda_0}^\lambda}\colon \widetilde{\mathcal{O}_{\lambda_0}}\to \widetilde{\mathcal{O}_\lambda}$ and $\widetilde{T_{\lambda}^{\lambda_0}}\colon \widetilde{\mathcal{O}_{\lambda}}\to \widetilde{\mathcal{O}_{\lambda_0}}$, respectively.

Using the argument in \cite{MR2005290}, we can prove the following properties.
Put $\widetilde{\theta} = \widetilde{T_\lambda^{\lambda_0}}\widetilde{T_{\lambda_0}^\lambda}$.
Set $W_{S_\lambda} = \Stab_W(\lambda)$ and let $w_\lambda\in W_{S_\lambda}$ be the longest element.
Then we can take $\widetilde{T_{\lambda_0}^\lambda}$ and $\widetilde{T_\lambda^{\lambda_0}}$ such that $\widetilde{\theta}$ is self-dual and there exists a natural transformation $\Id\langle \ell(w_\lambda)\rangle \to \widetilde{\theta}$ and $\widetilde{\theta}\to \Id\langle -\ell(w_\lambda)\rangle$.
Define a subset $W(S_\lambda)$ of $W$ by $W(S_\lambda) = \{x\in W\mid \text{$xs > x$ for all $s\in S_\lambda$}\}$.
Then for $x\in W(S_\lambda)$, $\widetilde{\theta}(\widetilde{M}^0(x\lambda_0))$ has a filtration $M_i$ such that $M_i/M_{i - 1}$ is isomorphic to $\bigoplus_{\ell(w) = i,w\in W_{S_\lambda}}\widetilde{M}^0(xw\lambda_0)\langle \ell(w_\lambda) - 2\ell(w)\rangle$.

\begin{proof}[Proof of Theorem~\ref{thm:translation functor, surjective, kernel}]
We prove \ref{enum:map by translation functor is surjective}.
By Theorem~\ref{thm:calculation of V}, for $y\in W(S_\lambda)$ and $w\in W_{S_\lambda}$, we have $V(x,yw)\subset V(x,y)\subset V(x,yw) + \sum_{s\in S_\lambda}Kv_s$.
It is easy to see that $v_s$ is in the kernel of $V_{\lambda_0}\to V_\lambda$.
Hence we may assume that $y \in W(S_\lambda)$.

It is sufficient to prove that 
\[
	\Ext^1(\widetilde{M}^0(x\lambda_0),\widetilde{M}^0(y\lambda_0)\langle 2\rangle)_0\to 
	\Ext^1(\widetilde{M}^0(x\lambda_0),\widetilde{\theta}(\widetilde{M}^0(y\lambda_0))\langle 2 - \ell(w_\lambda)\rangle)_0
\]
is surjective.
Let $M$ be the cokernel of 
\[
	\widetilde{M}^0(y\lambda_0)\langle 2\rangle\to \widetilde{\theta}(\widetilde{M}^0(y\lambda_0))\langle 2 - \ell(w_\lambda)\rangle.
\]
Then it is sufficient to prove that $\Ext^1(\widetilde{M}^0(x\lambda_0),M)_0 = 0$.
As we mentioned above, $M$ has a filtration $\{M'_i\}_{i\ge 1}$ such that $M'_i/M'_{i - 1}\simeq \bigoplus_{\ell(w) = i,w\in W_\lambda}\widetilde{M}^0(yw\lambda_0)\langle 2-2\ell(w)\rangle$.
By Lemma~\ref{lem:vanish of ext1, by grading}, we have $\Ext^1(\widetilde{M}^0(x\lambda_0),\widetilde{M}^0(yw\lambda_0)\langle 2 - \ell(w)\rangle)_0 = 0$ if $\ell(w) > 0$.
Hence we get \ref{enum:map by translation functor is surjective}.

We have \ref{enum:kernel of the map by translation functor} from \ref{enum:map by translation functor is surjective} and \ref{lem:preliminary of ext1:dim of V} of Lemma~\ref{lem:preliminary of ext1}.
\end{proof}

\def\cprime{$'$}

\clearpage

\section{Erratum to ``First extension groups of Verma modules and $R$-polynomials''}
Theorem 1.2 in \cite{MR3346063} (Theorem~\ref{thm:Main Theorem:Gabber-Joseph conjecture} and \ref{thm:Gabber-Joseph conjecture} in the main body) is false.
We can only prove the inequality, namely we only have the following theorem.
\begin{Theorem}
Assume that $\lambda$ is regular.
Then $\dim V_\lambda(x,y)$ is less than or equal to the coefficient of $q$ in $(-1)^{\ell(y) - \ell(x) - 1}R_{y,x}(q)$.
\end{Theorem}

When $x = w_{0}$, we have the equality.

\begin{Theorem}\label{thm:Gabber-Joseph conjecture, longest case}
Assume that $\lambda$ is regular.
Then $\dim V_\lambda(w_{0},x)$ is equal to the coefficient of $q$ in $(-1)^{\ell(x) - \ell(w_{0}) - 1}R_{x,w_{0}}(q)$.
\end{Theorem}

Both theorems follow from the proof of Theorem~\ref{thm:Gabber-Joseph conjecture}.
In the proof of Theorem~\ref{thm:Gabber-Joseph conjecture}, we claimed that Theorem~\ref{thm:Gabber-Joseph conjecture} follows from Theorem~\ref{thm:Gabber-Joseph conjecture, longest case}.
However this argument is not correct.

Here is a counterexample of Theorem~\ref{thm:Gabber-Joseph conjecture}.
We use the notation in the proof of Theorem~\ref{thm:Gabber-Joseph conjecture}.
Let $\mathfrak{g}$ be the simple Lie algebra of type $B_{3}$ and we use the standard notation of the root system.
In particular, the set of simple roots is $\{e_{1} - e_{2},e_{2} - e_{3},e_{3}\}$.
Let $s_{1} = s_{e_{1} - e_{2}},s_{2} = s_{e_{2} - e_{3}},s_{3} = s_{e_{3}}$ be simple reflections.
Put $x = s_{2}s_{3}s_{2}s_{1}s_{2}s_{3}$ and $y = s_{3}s_{2}$.
Then using the formula in the proof of Theorem~\ref{thm:Gabber-Joseph conjecture}, we have
\begin{align*}
r_{s_{2}s_{3}s_{2}s_{1}s_{2}s_{3},s_{3}s_{2}} & = r_{s_{2}s_{3}s_{2}s_{1}s_{2},s_{3}s_{2}} + 1\\
& = r_{s_{2}s_{3}s_{2}s_{1},s_{3}} + 1\\
& = r_{s_{2}s_{3}s_{2},s_{3}} + 2\\
& = r_{s_{2}s_{3},s_{3}} + 3\\
& = r_{s_{2},e} + 3\\
& = r_{e,e} + 4 = 4.
\end{align*}
However $n_{x,y}\le \dim \mathfrak{h} = 3$.
Hence $n_{x,y} < r_{x,y}$.

\subsection*{Acknowledgments}
Kevin Carlin pointed out that Theorem~\ref{thm:Main Theorem:Gabber-Joseph conjecture} implies the Gabber-Joseph conjecture for $\Ext^{1}$.
In particular, by Lemma~\ref{lem:preliminary of ext1}, for any $x\ge y$, we have $\Hom(M(x\lambda),M(\lambda)/M(y\lambda)) = 0$.
However, one can find a counterexample of this claim using \cite{MR3209798}.
The author thanks to Kevin Carlin for this comment.
The author also thanks to Hisayosi Matumoto who explained an example of $\Hom(M(x\lambda),M(\lambda)/M(y\lambda)) \ne 0$.

\end{document}